\documentclass{article}

\newtheorem{theorem}{Theorem}[section]
\newtheorem{lemma}[theorem]{Lemma}
\newtheorem{corollary}[theorem]{Corollary}
\newtheorem{definition}[theorem]{Definition}

\newcommand\nn{{\{0,1\}^n}}

\newcommand\lz{[\mkern-2.5mu[}
\newcommand\rz{]\mkern-2.5mu]}
\newcommand\lb{\langle\mkern-2.5mu\langle}
\newcommand\rb{\rangle\mkern-2.5mu\rangle}

\newcommand\ba{{\cal A}}
\newcommand\bb{{\cal B}}

\newcommand\qed{\begin{flushright} {\bf q.e.d.} \end{flushright} }
\newcommand\prf{\noindent {\bf Proof :}  }

\begin{document}

\title{A saturation property of structures obtained
by forcing with a compact family of random variables}

\author{Jan Kraj\'{\i}\v{c}ek\thanks{
Supported
in part by grant IAA100190902.
Also partially affiliated with the Institute of Mathematics of the Academy of Sciences.}}

\date{Faculty of Mathematics and Physics\\
Charles University in Prague\\
Sokolovsk\' a 83, Prague 8, CZ - 186 75\\
The Czech Republic\\
{\tt krajicek@karlin.mff.cuni.cz}
}

\maketitle

\begin{abstract}
A method for constructing Boolean-valued models 
of some fragments of arithmetic was developed in \cite{k2},
with the intended applications in bounded arithmetic and
proof complexity.
Such a model is formed by a family of random variables 
defined on a pseudo-finite sample space. 
We show that under a fairly 
natural condition on the family (called compactness in \cite{k2})
the resulting structure has a property that is naturally interpreted
as saturation for existential types.
We also give an example showing that this cannot be
extended to universal types.

\end{abstract}

Let $K$ be a Boolean-valued $L$-structure. That is, each 
sentence $A$ in the language $L(K)$, $L$ augmented by constants
for all elements of $K$, is assigned a truth-value
$\lz A \rz$ in a complete Boolean algebra $\bb$. These values
commute with propositional connectives (after Boole\cite{Boo})
and satisfy
$$
\lz \exists x A(x) \rz\ =\ 
\bigvee_{u \in K} \lz A(u) \rz\ \ \ \mbox {and }\ \ \ 
\lz \forall x A(x) \rz\ =\ 
\bigwedge_{u \in K} \lz A(u)\rz\ 
$$
(after Rasiowa-Sikorski\cite{RS1}).

Let ${\mathbf p}$ be a set of formulas 
in variables $x = x_1,\dots, x_n$. In the classical case
the set ${\mathbf p}$ is 
an {\bf $n$-type} over a structure if it is {\bf finitely satisfied}, i.e.
$$
\bigwedge_{\Phi \subseteq_{f} {\mathbf p}} \exists x 
\bigwedge_{A \in \Phi} A(x)
$$ 
where $\Phi$ runs over finite subsets of $\mathbf p$, and it is {\bf realized}
in the structure if
$$
\exists x \bigwedge_{A \in {\mathbf p}} A(x)
$$
holds there. In the context of Boolean-valued structures
this is naturally transcribed as the following
condition on truth values:
$$
({\sf Satur})\ \ \ \ \ \ \ \ \ \ \ \ 
\bigwedge_{\Phi \subseteq_{f} {\mathbf p}} \lz \exists x 
\bigwedge_{A \in \Phi} A(x)\rz \ \ \le\ \  
\bigwedge_{A \in {\mathbf p}} \lz A(u)\rz\ \ \ 
\ \ \ \ \ \ \ \ 
$$
for some $u \in K^n$; we shall say that such $u$ 
{\bf realizes ({\sf Satur}) for ${\mathbf p}$} 
in the structure. Note that in that case ({\sf Satur})
is actually an equality as the left-hand side always majorizes the
right-hand side.

In this paper we show that a certain class of
Boolean-valued structures constructed via forcing with random
variables (recalled in Sections \ref{1} and \ref{2}) is 
saturated in this sense for sets of existential formulas. 
This is done first
for the special case of sets ${\mathbf p}$ consisting
of open formulas in Section \ref{3} to display clearly the
idea, and for the existential case in Section \ref{4}.
In Section \ref{5} we give an example of a structure 
from the same class
that is not saturated for sets of universal formulas.
The paper is concluded by a brief explanation in
Section \ref{6} how is open
saturation of the structures considered
potentially relevant to proof complexity.

Background on model theory can be found in \cite{Mar},
further material (and details) on forcing with random variables
in \cite{k2}.

\section{Forcing with random variables set-up}
\label{1}

In this section we shall 
briefly recall the construction of Boolean valued structures
by forcing with random variables from \cite{k2}. The intended target structures
are  models of arithmetic with a special emphasis on bounded arithmetic.
This is motivated by a close relation of bounded arithmetic to
proof complexity but we shall not review this topic
here (an interested reader can consult \cite{kniha,k2}).

The structures are built from a family of random variables
on a pseudo-finite sample space. Let $\cal M$ be a non-standard 
$\aleph_1$-saturated model of true arithmetic in some language $L$
containing the language of Peano arithmetic
and having a canonical interpretation in the standard model ${\mathbf N}$. 
In \cite{k2} we took a language having symbols for all relations and
functions on ${\mathbf N}$ but here it is natural to 
consider countable $L$. In fact, we shall assume that $L$ is
definable in Peano arithmetic (this includes finite and recursive languages). 

Let $\Omega \in {\cal M}$
be an infinite set; as it is an element of the model it is $\cal M$-finite. 
Let $F \subseteq {\cal M}$ be any family 
of functions $\alpha : \Omega \rightarrow {\cal M}$. We call elements of
$F$ random variables. It is not assumed that $F$ is definable in the ambient
model ${\cal M}$.

Let $\ba$ be the Boolean algebra
of $\cal M$-definable subsets of $\Omega$ and let $\bb$ 
be its quotient by the ideal ${\cal I}$ of sets
of an infinitesimal counting measure.
Using the idea of Loeb's measure, the $\aleph_1$-saturation
of ${\cal M}$ and some measure theory
it was shown in \cite{k2} that $\bb$ is a complete Boolean algebra.

The counting measure on $\ba$ induces {\bf a strict measure $\mu$}
(in the ordinary sense with values in ${\mathbf R}$)
on $\bb$. The measure defines a metric on $\bb$: the distance of
two elements is the measure of their symmetric difference.

For any $k$-ary function symbol $f$ from $L$ and any
$\alpha_1, \dots , \alpha_k \in F$ define
the function $f(\alpha_1, \dots , \alpha_k) : \Omega \rightarrow {\cal M}$
by
$$
f(\alpha_1, \dots , \alpha_k)(\omega)\ :=\ 
f(\alpha_1(\omega), \dots , \alpha_k(\omega)), \ \mbox{ for } 
\omega \in \Omega\ .
$$
If this function is also always in $F$ we say that $F$ is {\bf $L$-closed}.

Any $L$-closed family $F$ is the universe of a Boolean-valued $L$-structure
$K(F)$ with $L(F)$-sentences having their
truth values in $\bb$ defined as follows.
Every atomic $L(F)$-sentence $A$ is naturally assigned a set 
$\lb A \rb$ from
$\ba$ consisting of those samples 
$\omega \in \Omega$ for which $A$ is true in ${\cal M}$.
The image of $\lb A \rb$ in $\cal B$, the quotient
$\lb A \rb / {\cal I}$, is denoted
$\lz A \rz$. Following Boole \cite{Boo} and Rasiowa-Sikorski \cite{RS1}
this determines the truth value $\lz A \rz \in {\cal B}$
for any $L(F)$-sentence $A$: $\lz \dots \rz$ commutes with
Boolean connectives and
$$
\lz \exists x A(x)\rz\ :=\ \bigvee_{\alpha\in F} \lz A(\alpha)\rz
\ \mbox {and }\ 
\lz \forall x A(x)\rz\ :=\ \bigwedge_{\alpha\in F} \lz A(\alpha)\rz\ .
$$
It holds that all logically valid sentences get the 
maximal truth value $1_{\bb}$. We say that a sentence is {\bf valid
in ${\cal M}$} if its truth value is $1_{\bb}$.

There are various generalizations of this basic set-up
considered in \cite{k2}. For example, 
the random variables from the family $F$
can be only partially defined on the sample space $\Omega$
(as long as their regions of undefinability have infinitesimal counting
measures) or the sample space may be equipped with some other
probability distribution than the uniform one.

We shall use 
one immediate consequence of the $\aleph_1$-saturation
of $\cal M$ and so we formulate it as a lemma.

\begin{lemma}
Let $\{a_k\}_{k \in {\mathbf N}}$ be any sequence of elements of ${\cal M}$.
Then there is an element $a^* \in {\cal M}$ that codes a sequence
$\{a^*_i\}_{i \le t}$ of some non-standard length $t \in {\cal M} \setminus {\mathbf N}$
such that $a^*_k = a_k$ for all $k \in {\mathbf N}$. 
\end{lemma}
Any such element $a^*$ is called {\bf a non-standard extension} of $\{a_k\}_{k \in {\mathbf N}}$.
We shall skip in future the ${\ }^*$ in the notation and denote a non-standard
extension simply as $\{a_i\}_{i \le t}$. 

\section{Compact families and witnessing of quantifiers}
\label{2}
In this section we shall recall the 
concept of a compact family $F$ from \cite[Chpt.3]{k2}
and some properties of $K(F)$ it implies.

\begin{definition} \label{2.1}
Let $F \subseteq {\cal M}$ be a family. 

\begin{enumerate}

\item 
$F$ is {\bf closed under definitions by
cases by open $L$-formulas} iff whenever $\alpha, \beta \in F$
and $B(x)$ is an open $L(F)$-formula with free variable $x$
then there is $\gamma \in F$ such that:

\[  \gamma(\omega) =
\left\{ \begin{array}{ll}
     \alpha(\omega)  &  \mbox{if $B(\alpha(\omega))$ holds}
 \\
     \beta(\omega)  &  \mbox{otherwise.}
                                      \end{array}
                              \right. \]

\item 

$F$ is {\bf compact}
iff there exists an $L$-formula $H(x,y)$ such that
for
$$
F_a := \{b \in {\cal M}\ |\ {\cal M} \models H(a,b)\}
$$
the following two properties hold:
\begin{itemize}

\item $\bigcap_{k \in {\mathbf N}} F_k = F$.

\item $F_k \supseteq F_{k+1}$, for all $k \in {\mathbf N}$.
\end{itemize}
\end{enumerate}
\end{definition}
Recall that the Overspill is the principle, a simple consequence of induction,
that if a property definable in ${\cal M}$ holds for all standard
numbers, it must hold actually for all elements up to some non-standard
element of ${\cal M}$ (cf. \cite[Appendix]{k2}).
The primary intended use of compact families is to allow the following 
type of reasoning.

Assume $\{\alpha_k\}_{k \in \mathbf N}$ is an arbitrary sequence
of elements of $F$ and that 
$\{\alpha_i\}_{i \le t}$ is its non-standard extension. 
Note that the conditions posed on sets $F_k$ in the definition of
compactness imply that for all standard $k$ it holds
that 
\begin{itemize}
\item $\forall j \le k \ \alpha_j \in F_k$, and
\item $\forall j < k \ F_j \supseteq F_{j+1}$.
\end{itemize}
By the Overspill in ${\cal M}$ this must hold also for some non-standard
$s \le t$. In particular, all $\alpha_j$ from the non-standard extension
with $j \le s$ are in $F_s \subseteq F$, and hence in $F$ too.

\begin{theorem} \label{2.2} {\cite[Thm.3.5.2]{k2}}

Let $F$ be an $L$-closed family that is
closed under definition by cases by open
$L$-formulas and compact.
Let $A$ be an $L(F)$-sentence of the form 
$$
\exists x_1 \forall y_1 \dots
\exists x_k \forall y_k B(x_1, y_1, \dots, x_k, y_k)
$$
with $B$ open.

Then there are $\alpha_1, \beta_1, \dots,
\alpha_k, \beta_k \in F$ such that for all $i = 1, \dots, k$:
 
$$
\lz \forall y_i \exists x_{i+1} \forall y_{i+1} \dots
\exists x_k \forall y_k 
B(\alpha_1, \beta_1, \dots,
\alpha_i, y_i,x_{i+1}, y_{i+1},\dots, x_k, y_k) \rz
= \lz A \rz\ 
$$ 
and
$$
\lz \exists x_{i+1} \forall y_{i+1} \dots
\exists x_k \forall y_k 
B(\alpha_1, \beta_1, \dots,
\alpha_i, \beta_i,x_{i+1}, y_{i+1},\dots, x_k, y_k) \rz
= \lz A \rz\ .
$$

\end{theorem}

\section{Saturation for sets of open formulas}
\label{3}

In this section we prove a special case of saturation
when all formulas in the set ${\mathbf p}$ are open. We will use
the following immediate corollary of Theorem \ref{2.2}.
To simplify the notation we consider, here as well as
in the next section, sets of formulas in one free
variable; the general case would be done in the same way.

\begin{corollary} \label{3.1}

Let $F$ be an $L$-closed family that is
closed under definition by cases by open
$L$-formulas and compact.
Let $A(x)$ be an open $L(F)$-formula with one free variable $x$.

Then there is $\alpha \in F$ such that:
$$
\lz \exists x A(x) \rz
= \lz A(\alpha) \rz\ =\ \lb A(\alpha)\rb/{\cal I}
\ . 
$$ 
\end{corollary}

\prf

The first equality follows from Theorem \ref{2.2}
and the second one follows by the definition of
$\lz A \rz$ as $A$ is open and taking the quotient by
${\cal I}$ commutes with Boolean connectives.

\qed

\medskip
In proving the next theorem we could restrict to the special
case when the left-hand side of ({\sf Satur}) has value
$1_{\bb}$ by taking a suitable quotient of $\bb$. However,
this would be done at the expense of having to show that
the inequality for the resulting new Boolean-valued structure 
could be pulled back to the original one.
We thus prefer not to make this simplification.

\begin{theorem} \label{3.2}

Let $L$ be a language definable in Peano arithmetic.
Let $F$ be an $L$-closed family that is
closed under definition by cases by open
$L$-formulas and compact.

Assume that ${\mathbf p}$ is a countable set of open 
$L(F)$-formulas in one free variable $x$.

Then there is an element $\alpha \in F$ that realizes 
the saturation inequality $({\sf Satur})$ for the set 
${\mathbf p}$ in $K(F)$.

\end{theorem}

\prf

Let $A'_1(x), A'_2(x), \dots$ enumerate ${\mathbf p}$ and define
$$
A_k(x) \ :=\ \bigwedge_{i \le k} A'_i(x)\ ,\ \mbox{ for $k \geq 1$}\ .
$$
Then we have the following

\medskip
\noindent
{\bf Claim 1:}{\em All implications
$$
A_{k+1}(x) \rightarrow A_k(x)
$$
are logically valid and 
if an element $\alpha \in F$ satisfies 
the following inequality
$$
\bigwedge_k \lz \exists x A_k(x)\rz\ \le\ 
\bigwedge_k \lz A_k(\alpha)\rz
$$
then $\alpha$ realizes the $({\sf Satur})$ inequality for ${\mathbf p}$
in $K(F)$.}

\medskip
\noindent
Note that the inequality is actually an equality in that case.

\medskip
\noindent
{\bf Claim 2:}{\em 
Assume $\alpha \in F$, $U \subseteq \Omega \wedge U \in {\cal M}$ and it holds
$$
\lb A_k(\alpha)\rb \supseteq U
$$
for all $k \in {\mathbf N}$, and
$$
\mu(\lz \exists x A_k(x)\rz)\ \searrow\ 
\mu(U/{\cal I})
$$
as standard $k \rightarrow \infty$.

Then $\alpha$ satisfies the inequality 
from Claim 1 and hence realizes ({\sf Satur}) for ${\mathbf p}$
in $K(F)$.}

\medskip
\noindent
The arrow $\searrow$ means that the sequences of reals
on the left-hand side is non-increasing and its limit is 
the right-hand side.
The claim follows as
$$
\lz \exists x A_k(x)\rz)\ \geq\ 
\lz A_k(\alpha)\rz)\ \geq\ 
U/{\cal I}\ .
$$

\bigskip

Now we are going to show that some $\alpha$ and $U$ satisfying the
hypothesis of Claim 2 do exist. This will prove the theorem.

By Corollary \ref{3.1} there are $\alpha_k \in F$ such that
$$
\lz \exists x A_k(x)\rz\ =\ \lz A_k(\alpha_k)\rz\ .
$$
Define $U_k := \lb A_k(\alpha_k)\rb$, so
$$
\lz \exists x A_k(x)\rz\ =\ U_k/{\cal I}
$$
and as $A_{\ell}$ logically implies $A_k$ for $l\geq k$, also
$$
U_k\ \supseteq \ U_{\ell}\ ,\ \mbox{ for }\ell \geq k\ .
$$

\bigskip

Consider the sequence $\{A_k, \alpha_k, U_k\}_{k \in {\mathbf N}}$ 
and let
$$
\{A_i, \alpha_i, U_i\}_{i \le t}\ ,\ \mbox{ $t$ non-standard } 
$$
be its non-standard extension provided by Lemma \ref{2.1}.

By the compactness of $F$ there is a definable family $\{F_a\}_a$
of sets such that $F = \bigcap_{k\ \in {\mathbf N}} F_k$.
Consider the following property, definable in ${\cal M}$
(the definability of $L$ is used here), of an element
$i \le t$.
It is the conjunction of seven conditions:
\begin{enumerate}

\item $A_i$ is an open $L(F_i)$-formula.

\item $A_i \rightarrow A_j$ is logically valid for all $j \le i$.

\item $\lb A_i(\alpha_i)\rb = U_i$.

\item $\Omega \supseteq U_1 \supseteq \dots \supseteq U_i$.

\item $\frac {|U_j|}{|\Omega|} - 
\frac {|U_i|}{|\Omega|} < 1/j$ for all $j \le i$.

\item $F_1 \supseteq \dots \supseteq F_i$.

\item $\alpha_j \in F_i$ for all $j \le i$.

\end{enumerate}
All seven conditions are valid for all standard $i$ 
perhaps with the exception of 5: but taking a suitable
subsequence of the
original sequence $\{A_k, \alpha_k, U_k\}_{k \in {\mathbf N}}$ 
arranges this condition too.
The first two items are included because in the argument below
we need to talk in $\cal M$
also about the satisfiability relation for formulas $A_i$ with 
a non-standard index.
By the Overspill then the property must be true for all $i \le s$ 
up to some non-standard $s \le t$.
We want to show that for such an $s$, $\alpha_s$ and $U_s$
satisfy the hypothesis of Claim 2.

\bigskip

By 6 and 7 $\alpha_s \in F_s \subseteq F$ and by 
1, 2 and 3
$$
\lb A_k(\alpha_s)\rb\ \supseteq\ \lb A_s(\alpha_s)\rb\ =\ U_s
$$
for all $k \in {\mathbf N}$. It remains to note that
$$
\mu(\lz \exists x (A_k(x)\rz)\ \searrow\
\mu(U_s/{\cal I})\ .
$$
follows by 3, 4 and 5.

\qed

In the argument we have used that for an open formula $A(x)$ 
the set $\lb A(\alpha)\rb$ is definable from $\alpha$
and satisfies $\lz A(\alpha)\rz = \lb A(\alpha)\rb/{\cal I}$.
This is not true for general formulas but for for existential formulas
one could add witnesses (in the sense of Corollary \ref{3.1})
for the values $\lz A_k(\alpha_k)\rz$ to the data 
in $\{A_k, \alpha_k, U_k\}_{k \in {\mathbf N}}$ and run an
analogous argument.

Assuming little bit more about the structure however,
we can derive the existential case directly from Theorem
\ref{3.2}.

\section{The existential case}
\label{4}

In this section we note that Theorem \ref{3.2} implies the
statement also for sets of existential formulas as long as
the underlying structure admits a pairing function. This is always
the case for structures of interest in \cite{k2} as they are
models of various bounded arithmetics.

\begin{definition}
Let $F$ be an $L$-closed family. We say that $K(F)$ {\bf has
pairing} if $L$ contains symbols $\pi_1(x), \pi_2(x)$
for two unary functions and a symbol $\langle x, y\rangle$ 
for a binary function, and the
universal closures of the following three formulas are valid
in $K(F)$:
\begin{itemize}

\item $\langle \pi_1(z), \pi_2(z)\rangle = z$.

\item $\pi_1(\langle x, y\rangle) = x$.

\item $\pi_2(\langle x, y\rangle) = y$.

\end{itemize}
\end{definition}

\begin{theorem} \label{4.2}

Let $L$ be a language definable in Peano arithmetic.
Let $F$ be an $L$-closed family that is
closed under definition by cases by open
$L$-formulas and compact. Assume that $K(F)$ has pairing.

Let ${\mathbf p}$ be a countable set of existential 
$L(F)$-formulas in one free variable $x$.
Then there is an element $\alpha \in F$ that realizes 
the saturation inequality $({\sf Satur})$ for the set 
${\mathbf p}$ in $K(F)$.

\end{theorem}

\prf

Using the pairing to replace several existential quantifiers
by one we may assume without a loss of generality
that the formulas in ${\mathbf p}$ have the form
$A(x) = \exists y B(x,y)$, with $B$ open. Enumerate
${\mathbf p}$ as $A_k(x) = \exists y B_k(x,y)$, $k \in {\mathbf N}$,
and define open formulas
$$
C_k(z)\ :=\ 
B_k(\pi_1(z), \pi_1(\pi_2^{(k)}(z)))\ 
$$
where $\pi_2^{(k)}$ abbreviates $k$-times iterated
$\pi_2$.
The following claim implies that if $\alpha \in F$ realizes
({\sf Satur}) for $\{C_k\ |\ k\in {\mathbf N}\}$ then
$\pi_1(\alpha)$ realizes it for ${\mathbf p}$.

\medskip
\noindent
{\bf Claim:} {\em
For any $k \in {\mathbf N}$ it holds that
$$
\lz \exists x \bigwedge_{i \le k} A_i(x)\rz\ =\ 
\lz \exists z \bigwedge_{i \le k} C_i(z)\rz\ .  
$$
} 

\medskip

Clearly the left-hand side majorizes the right-hand side.
For the opposite direction apply Corollary \ref{3.1}
to get $\gamma, \beta_1, \dots, \beta_k \in F$
such that
$$
\lz \exists x \bigwedge_{i \le k} A_i(x)\rz\ =\ 
\bigwedge_{i \le k} \lz B_i(\gamma, \beta_i)\rz\ 
$$
and define
$$
\alpha\ :=\ 
\langle \gamma, \langle \beta_1, \langle \beta_2 , \dots,
\langle \beta_k, \gamma
\rangle \dots \rangle
$$
(the second $\gamma$ could be replaced by any element of $F$.)
It is easy to see that the substitution $z := \alpha$ gives 
to the right-hand side in the Claim a value that equals to 
$\lz \exists x \bigwedge_{i \le k} A_i(x)\rz$.

\qed

\section{The failure of the universal case}
\label{5}

In this section we show that Theorem \ref{4.2}
cannot be generally strengthened to sets of universal
formulas. Take for $L$ the language of Peano arithmetic
together with the inequality sign $\le$ and with some function 
symbols for pairing and its projections, and having also a 
unary function symbol $|x|$ for the bit-length of 
number $x$.

In ${\cal M}$ we shall identify numbers with the binary
strings consisting of their bits.
Let the sample space $\Omega$ be simply $\nn$ for some
non-standard $n \in {\cal M}$ and let the family $F$ 
consist of all functions on $\Omega$ computed by
circuits with $n$ inputs and arbitrarily many outputs
but of the size bounded above by all terms
$2^{n^{1/k}}$, for all $k\in {\mathbf N}$.
In other words, these are functions computed in 
sub-exponential non-uniform time.
It is easy to see that $F$ is compact: take for $F_k$
the functions on $\Omega$ computed by
circuits of size bounded above by $2^{n^{1/k}}$.

Let $id_{\Omega}$ be the identity function on $\Omega$.
Consider the following universal $L$-formulas:
$$
A_k(x)\ :=\ 
|x|^k \le |id_{\Omega}| \wedge \forall y (|y| \neq x)\ .
$$
(The expression $|x|^k$ is just an abbreviation for 
the term $|x|\cdot \dots \cdot |x|$, $|x|$ occurring 
$k$-times.)

The value of $|id_{\Omega}|$ is on all samples $n$.
Take for $x$ any constant function $\alpha$ outputting 
a fixed string of bit-length $n^{1/k}$ and of
value $2^{n^{1/k}}$. Then clearly
$\lz |\alpha|^k \le |id_{\Omega}|\lz = 1_{\bb}$
but also $\lz \forall y (|y| \neq \alpha)\rz = 1_{\bb}$.
This is because
no function from $F$ can output a string of length
$2^{n^{1/k}}$. This implies that 
the left-hand side of ({\sf Satur})
for formulas $A_k(x)$ has the truth value $1_{\bb}$.

On the other hand, assume that $\alpha$ makes all
sentences $|\alpha|^k \le |id_{\Omega}|$ valid.
Hence the sets $U_k \subseteq \Omega$ defined
by
$$
\omega \in U_k \ \mbox { iff }\ 
|\alpha(\omega)| \le n^{1/k}
$$
have all counting measures infinitesimally close to
$1$. The sequence of $U_k$ thus satisfies
$$
U_1 \supseteq \dots \supseteq U_k
$$
and 
$$
1\ -\ \frac{|U_i|}{|\Omega|} \ <\ 1/k\ ,\ \mbox{ for all $i \le k$}\ .
$$
Taking its non-standard extension $\{U_i\}_{i \le t}$
and applying the Overspill analogously 
as before yields a non-standard
$s \le t$ for which
the counting measure of 
$U_s = \lb |\alpha|^s \le |id_{\Omega}|\rb$ is
infinitesimally close to $1$, i.e. 
$\lz |\alpha|^s \le |id_{\Omega}|\rz = 1_{\bb}$ holds. 

But then the value of $\alpha$
on each sample from $U_s$ is at most $2^{n^{1/s}}$ and
so there is a $\beta \in F$
that outputs on each sample $\omega \in U_s$ 
a string of bit-length
$\alpha(\omega) \le 2^{n^{1/s}}$. 
For such $\beta$ however, 
$\lz |\beta| \neq \alpha \rz = 0_{\bb}$. 
This argument proves the following statement.

\begin{theorem}

There is a finite language $L$ and an $L$-closed family $F$ that is
closed under definition by cases by open
$L$-formulas, compact and such that $K(F)$ has pairing,
but for which there 
is a countable set of universal 
$L(F)$-formulas in one free variable $x$ for which no
element of $K(F)$ realizes 
the saturation inequality $({\sf Satur})$.

\end{theorem}

\section{A concluding remark}
\label{6}

In the abstract we have alluded to intended applications
of forcing with random variables in proof complexity.
To illustrate - in rather abstract terms - how open saturation
can be useful consider the following situation.

Let $T$ be a theory in a countable language $L$ and let
$\varphi$ be an $L$-definable 3CNF propositional 
formula $\varphi$. In a typical example $T$ may be an extension of 
a bounded arithmetic theory by an open diagram of some $L$-structure.
The task would be to construct a model of $T$ containing 
a truth assignment satisfying $\varphi$. That is, the model
should satisfy a sentence of the form 
$$
\exists x \forall y A(x,y)
$$
where $A(x,y)$ is an open formula
formalizing that the assignment $x$ satisfies the $y$-th 
3-clause of $\varphi$. 

If such a model can be found with $\alpha$ satisfying $\forall y A(\alpha,y)$, 
it often actually suffices to work further only
with a substructure of the model generated by $\alpha$. 
But in order for this substructure to
satisfy $\forall y A(\alpha,y)$ it is not necessary that the original model
satisfies $\forall y A(\alpha, y)$; it suffices 
that $\alpha$ realizes in the model the open type consisting 
of countably many formulas
$$
A(x,t(x))
$$
with $t(x)$ ranging over all (countably many) terms with the only  
free variable $x$.

Open saturation allows to simplify the task to construct such a model
and to consider only how to realize in it finite subsets of the type.

\end{document}